\documentclass[11pt]{article}

\usepackage[all]{xy}
\usepackage{amsmath}
\usepackage{amsfonts}
\usepackage{amssymb}
\usepackage{amscd}
\usepackage{amsthm}
\usepackage{latexsym}
\usepackage{amsbsy}
\usepackage{epsfig, setspace}
\usepackage{subfigure}
\usepackage{graphicx}

\newtheorem{theorem}{Theorem}[section]
\newtheorem{remark}[theorem]{Remark} 

\newtheorem{lemma}[theorem]{Lemma}

\newtheorem{example}[theorem]{Example}

\newcommand{\rnc}[2]{\renewcommand{#1}{#2}}
\rnc{\theequation}{\thesection.\arabic{equation}}
\rnc{\[}{\begin{equation}}
\rnc{\]}{\end{equation}}

\begin{document}
\title{Enveloping algebras of some quantum Lie algebras} 

\author{Arash Pourkia\\ {\small Huron University College, Western University,}\\ {\small Department of Mathematics, Western University,}\\ {\small London, Ontario, Canada}} 

\date{August 2014}

\maketitle

\begin{abstract}
We define a family of Hopf algebra objects, $H$, in the braided category of $\mathbb{Z}_n$-modules (known as anyonic vector spaces), for which the property $\psi^2_{H\otimes H}=id_{H\otimes H}$ holds. We will show that these anyonic Hopf algebras are, in fact, the enveloping (Hopf) algebras of particular quantum Lie algebras, also with the property $\psi^2=id$. Then we compute the braided periodic Hopf cyclic cohomology of these Hopf algebras. For that, we will show the following fact: analogous to the non-super and the super case, the well known relation between the periodic Hopf cyclic cohomology of an enveloping (super) algebra and the (super) Lie algebra homology also holds for these particular quantum Lie algebras, in the category of anyonic vector spaces.
\end{abstract}

\section{Introduction} 
In \cite{pl4}, we proved that in the non-symmetric category of anyonic vector spaces there are families of examples in which, for any two objects $A$ and $B$, we have $\psi_{B\otimes A}\psi_{A\otimes B}=id_{A\otimes B}$. The motivation for this came from the fact that we showed that the main theorems of braided Hopf cyclic cohomology in \cite{pl1hha}, stated for symmetric monoidal categories, are still valid in non-symmetric categories as long as we have $\psi^2_{H \otimes H}=id_{H \otimes H}$ (or for short $\psi^2=id$) only for the object of  interest, $H$. We should recall that the braided Hopf cyclic theory in \cite{pl1hha,pl4} is a generalization of the Connes and Moscovici Hopf cyclic cohomology in \cite{cm2,cm3,cm4}, and of the more general case of Hopf cyclic cohomology with coefficients in \cite{hkrs1,hkrs2}, to the context of abelian braided monoidal categories. A short survey on the subject can be found in   \cite{pl3}.\newline 
 
This paper is organized as follows. In Section \ref{anyonhpfalgwtheproppsi2=id} we define a family of Hopf algebra objects, $H$, in the non-symmetric category of anyonic vector spaces, for which $\psi^2_{H \otimes H}=id_{H \otimes H}$. 
In Section \ref{aquantliealganditsenvalg} we construct a family of particular quantum Lie algebras, $\mathfrak{g}$, in the sense of \cite{maj-gom2003} (also with the property $\psi^2=id$), in the same category of anyonic vector spaces. Then, we will show that their universal enveloping algebras, $U(\mathfrak{g})$,  are the anyonic Hopf algebras introduced in Section \ref{anyonhpfalgwtheproppsi2=id}. 

We should mention that in \cite{maj-gom2003}, they have shown, among other things, that the notion of quantum Lie algebras, in the sense of \cite{wor89}, is closely related to braided Lie algebras, in the sense of  \cite{maj94-quanandbrddliealg94}. To do this, they have provided an axiomatic definition of (left) quantum Lie algebras. In this paper, we follow \cite{maj-gom2003} for the definition of a quantum Lie algebra, and \cite{maja94-algandhopfalginbrddcats, maj94-quanandbrddliealg94,maj95-book} for other definitions and notions we need. 

In Section \ref{hccforouranyonhofalg}  we present a lemma to compute the (braided) periodic Hopf cyclic cohomology of the anyonic Hopf algebras of Section \ref{anyonhpfalgwtheproppsi2=id}. It states that, analogous to the non-super \cite{cm2,cm3,cm4,crn} and the super case \cite{pl1hha}, the following well known relation between the periodic Hopf cyclic cohomology of an enveloping algebra, $HP_{(\delta,1)}^*(U(\mathfrak{g}))$, and the Lie algebra homologies, $H_i (\mathfrak{g} ; \mathbb{C}_ \delta)$, 
\[ HP_{(\delta,1)}^*(U(\mathfrak{g}))=\bigoplus_{i=*(mod~2)} H_i (\mathfrak{g} ; \mathbb{C}_ \delta), \,\, *=0,\,1, \,\,  (\delta \, {\text {is a character for}} \,\mathfrak{g}),   \nonumber\] 
holds for those quantum Lie algebras of Section \ref{aquantliealganditsenvalg}.   \newline

\thanks{The author would like to thank Masoud Khalkhali for his continuous encouragement and support.}

\section{Anyonic Hopf algebras with the property $\psi^2=id$} \label{anyonhpfalgwtheproppsi2=id}
In \cite{pl4}, we showed that in some of the non-symmetric categories of anyonic vector spaces, there are families of objects in which for any two objects $A$ and $B$ we have $\psi_{B\otimes A}\psi_{A\otimes B}=id_{A\otimes B}$. In this section, after a quick review, we provide examples of Hopf algebra objects within one of those families. In the next sections we set out to compute the Hopf cyclic cohomology of these Hopf algebras using the braided Hopf cyclic theory in \cite{pl1hha, pl4}.  \newline

Let $\mathbb{C} \mathbb{Z}_n$ be the  group (Hopf) algebra of the finite cyclic group, $\mathbb{Z}_n$, of order $n$. Let $R=R_1\otimes R_2$ be the nontrivial quasitriangular structure for $\mathbb{C} \mathbb{Z}_n$ defined by \cite{maja94-algandhopfalginbrddcats,maj95-book},
\[R=(1/n)\sum_{a,b=0}^{n-1} e^{\frac{(-2 \pi iab)}{n}} g^a \otimes g^b,  \label{Rforczn} \nonumber\] where $g$ is the generator of $\mathbb{Z}_n$. The category of all left $\mathbb{C} \mathbb{Z}_n$-modules, $\mathcal{C}$, is known as the category of anyonic vector spaces. The objects of $\mathcal{C}$ are of the form $V=\bigoplus_{i=0}^{n-1}\, V_i$. They are $\mathbb{Z}_n$-graded representations of $\mathbb{C} \mathbb{Z}_n$ and the action of $\mathbb{Z}_n$ on $V$ is given by, 
\[g\rhd v= e^\frac{2\pi i |v|}{n} v, \label{actofgonv} \nonumber\] 
where $|v|=k$ is the degree of the homogeneous elements $v$ in $V_k$. The morphisms of $\mathcal{C}$ are linear maps that preserve the grading. The braiding map in $\mathcal{C}$ is, 
\[\psi_{V \otimes W} (v \otimes w) = e^\frac{2\pi i |v||w|}{n}\,\, w\otimes v, \label{anyonbraid} \]
where $|v|$ and $|w|$ are the degrees of homogeneous elements $v$ and $w$ in objects $V$ and $W$, respectively. \newline

This category is the symmetric category of super vector spaces for $n=2$, but is not symmetric when $n>2$, \cite{maja94-algandhopfalginbrddcats,maj95-book}. However, as we proved in \cite{pl4}, for many values of $n>2$, one can always find objects $A$ and $B$ in $\mathcal{C}$ such that ``behave symmetric", i.e., $\psi_{B \otimes A}\psi_{A \otimes B}=id_{A \otimes B}$. The following is a family of such examples. 

\begin{example} \label{n=2m2}
Let $\mathcal{C}$ be the category of left $\mathbb{C} \mathbb{Z}_{n}$-modules, where $n=2m^2$ for some integer $m>1$. Let $A=\bigoplus_{i=0}^{n-1}\, A_i$ and $B=\bigoplus_{i=0}^{n-1}\, B_i$ be  objects in $\mathcal{C}$ which are focused only in degrees, $km$, for integers $k\geq 0$. By $A$ being focused only in these degrees we mean, $A_i = 0$ when $i\neq km$, for integers $k\geq 0$. For any two homogeneous elements $x$ in $A$ and $y$ in $B$, with $|x|=km$ and $|y|=lm$,\, $k,l \geq 0$, by Formula \eqref{anyonbraid}, we have, 
\[\psi_{A \otimes B} (x \otimes y) = e^{\frac{2\pi i |x||y|}{2m^2}} (y \otimes x)= (-1)^{\frac{|x||y|}{m^2}} (y \otimes x) = (-1)^{kl} (y \otimes x) . \label{siform2almostsuper} \]
Thus, it is clear that, 
\[ \psi_{B \otimes A}\psi_{A \otimes B}=id_{A \otimes B} \label{siBAsiAB=idAB}\]
 
An example of this case, as we will use in this paper, is when $n=18$ and $A=\bigoplus_{i=0}^{17}\, A_i$ and $B=\bigoplus_{i=0}^{17}\, B_i$ are objects focused only in degrees $0$, $3$, $6$, $9$, $12$ and $15$. 
\end{example}

\begin{remark} \label{remarkonweaklybrddobjects}
Formula \eqref{siform2almostsuper} shows that even though the objects in Example \ref{n=2m2}, $A$ and $B$, live in a truly braided category, their braid maps $\psi_{A \otimes B}$ behave almost like the braid maps for objects in the category of super vector spaces, by replacing $|x|$ and $|y|$ in the super case with $\frac{|x|}{m}$ and $\frac{|y|}{m}$.  
\end{remark}

Now we define a family of Hopf algebra objects, within the objects of Example \ref{n=2m2}. For such an anyonic Hopf algebra \cite{maja94-algandhopfalginbrddcats,maj95-book}, $H$, Formula \eqref{siBAsiAB=idAB} translates to $\psi^2_{H \otimes H}=id_{H \otimes H}$, as desired. 

\begin{example} \label{anyhopfalginz2m2moregeneral}
Let $\mathcal{C}$ be the category of left $\mathbb{C} \mathbb{Z}_{n}$-modules, where $n=2m^2$ for some integer $m>1$. We define a Hopf algebra object, $H=\bigoplus_{i=0}^{n-1}\, H_i$, in $\mathcal{C}$, focused in degrees $km$, for integers $k\geq 0$, as follows. As an algebra, we define $H$ to be generated by elements $x_i$ in degrees $im$ (i.e., $|x_i|=im$), for $i\geq 0$, subject to no relations (except, of course, $x_i^0=1$). We define the rest of the Hopf algebra structure on $H$ by, 
\[ \Delta(x_i)=1\otimes x_i + x_i\otimes 1, \quad \epsilon(x_i)=0, \quad S(x_i)=-x_i,\nonumber\]
$\Delta$ and $\epsilon$ extended as (braided) algebra maps, i.e., (using formula \eqref{siform2almostsuper}), for any homogeneous elements $v$ and $w$ in $H$, 
\[\Delta(vw)=(-1)^{\frac{|v^{(2)}||w^{(1)}|}{m^2}}v^{(1)}w^{(1)} \otimes v^{(2)}w^{(2)} \nonumber , \quad \epsilon(vw)=\epsilon(v)\epsilon(w)\]  and $S$ as an anti-commutative algebra map, i.e., 
\[S(vw)=(-1)^{\frac{|v||w|}{m^2}}S(w)S(v) .\nonumber  \] 
 This means, for example, if $|x_i|=im$ and $|x_j|=jm$ are both odd then we have, 
\[ \Delta(x_ix_j)=1\otimes x_ix_j -x_j\otimes x_i + x_i\otimes x_j+ x_ix_j\otimes 1, \quad \epsilon(x_ix_j)=0, \quad S(x_ix_j)= -x_jx_i,\nonumber\]
and if one of $|x_i|$ or $|x_j|$ is even, 
\[ \Delta(x_ix_j)=1\otimes x_ix_j +x_j\otimes x_i + x_i\otimes x_j+ x_ix_j\otimes 1, \quad \epsilon(x_ix_j)=0, \quad S(x_ix_j)= x_jx_i.\nonumber\]
\end{example}

\begin{remark}
In fact, in Example \ref{anyhopfalginz2m2moregeneral}, there is even no need to restrict ourselves to one generator for each degree. This means we can have, for example, a set of generators, $\{x_{i1}, x_{i2}, x_{i3}, \cdots \}$, in each degree $im$.
 \end{remark}

The structure defined in Example \ref{anyhopfalginz2m2moregeneral} is good for showing that, in general, we have a source of examples of anyonic Hopf algebras with the property $\psi^2=id$. However, it is too general. Below, we provide some simpler examples by working with an odd integer $m$, limiting the number of generators, and by imposing some relations into the algebra structure. They are easier, for our purposes, to work with and we will use them in the next sections. In particular, we consider the case $n=18$, i.e., $m=3$, very often.  

\begin{example} \label{anyhopfalginz2m2lessgeneral}
Let $\mathcal{C}$ be the category of left $\mathbb{C} \mathbb{Z}_{n}$-modules, where $n=2m^2$ for some odd integer $m>1$. We define a Hopf algebra object $H=\bigoplus_{i=0}^{n-1}\, H_i$, in $\mathcal{C}$, as follows. As an algebra, we define $H$ to be generated by elements $x_i$ in degrees $im$, for $i=1, 3, 5, 7, \cdots, 2m-1$, subject to the following relations. For $i,j=1, 3, 5, 7, \cdots, 2m-1$,
\[x_i^0=1, \qquad x_ix_j=-x_jx_i. \nonumber \]
We define the rest of Hopf algebra structure on $H$ by, 
\[ \Delta(x_i)=1\otimes x_i + x_i\otimes 1, \quad \epsilon(x_i)=0, \quad S(x_i)=-x_i,\nonumber\]
$\Delta$ and $\epsilon$ to be extended as (braided) algebra maps and $S$ as an anti-commutative algebra map. 
\end{example}

In particular if $n=18$, we have, 

\begin{example} \label{anyhopfalginz18lessgeneral}
Let $\mathcal{C}$ be the category of left $\mathbb{C} \mathbb{Z}_{18}$-modules. We define a Hopf algebra object $H=\bigoplus_{i=0}^{17}\, H_i$, in $\mathcal{C}$, as follows. For the algebra structure of $H$, we define $H$ to be generated by elements $x$, $y$, and $z$, of degrees $3$, $9$, and $15$, respectively, subject to the relations,
\[x^0=y^0=z^0=1, \quad xy=-yx, \quad xz=-zx, \quad yz=-zy. \nonumber \]
As regards the rest of the Hopf algebra structure on $H$, for $v=x,y,z$,  
\[ \Delta(v)=1\otimes v + v\otimes 1, \quad \epsilon(v)=0, \quad S(v)=-v,\nonumber\]
$\Delta$ and $\epsilon$ extended as (braided) algebra maps and $S$ as an anti-commutative algebra map. 
\end{example}

The following examples are yet simpler versions of the above ones, obtained by adding the relation $x_i^2=0$. 

\begin{example} \label{anyhopfalginz2m2}
Let $\mathcal{C}$ be the category of left $\mathbb{C} \mathbb{Z}_{n}$-modules, where $n=2m^2$ for some odd integer $m>1$. We define a Hopf algebra object, $H=\bigoplus_{i=0}^{n-1}\, H_i$, in $\mathcal{C}$, as follows. As an algebra, we define $H$ to be generated by elements $x_i$ in degrees $im$, for $i=1, 3, 5, 7, \cdots, 2m-1$, subject to the relations, for $i,j=1, 3, 5, 7, \cdots, 2m-1$,
\[x_i^2=0, \qquad x_i^0=1, \qquad x_ix_j=-x_jx_i. \nonumber \]
The rest of Hopf algebra structure on $H$ is defined by, 
\[ \Delta(x_i)=1\otimes x_i + x_i\otimes 1, \quad \epsilon(x_i)=0, \quad S(x_i)=-x_i,\nonumber\]
$\Delta$ and $\epsilon$ extended as (braided) algebra maps and $S$ as an anti-commutative algebra map. 
\end{example}

In particular if $n=18$, we have, 

\begin{example} \label{anyhopfalginz18}
The Hopf algebra object, $H=\bigoplus_{i=0}^{17}\, H_i$, in $\mathcal{C}$, the category of left $\mathbb{C} \mathbb{Z}_{18}$-modules, is defined as follows. As an algebra, $H$ is generated by elements $x$, $y$, and $z$, of degrees $3$, $9$, and $15$, respectively, subject to the relations,
\[x^2=y^2=z^2=0, \,\,\, x^0=y^0=z^0=1, \,\,\, xy=-yx, \,\,\, xz=-zx, \,\,\, yz=-zy. \nonumber \]
Regarding the rest of Hopf algebra structure on $H$, for $v=x,y,z$,  
\[ \Delta(v)=1\otimes v + v\otimes 1, \quad \epsilon(v)=0, \quad S(v)=-v,\nonumber\]
$\Delta$ and $\epsilon$ extended as (braided) algebra maps and $S$ as an anti-commutative algebra map. 
\end{example}

\section{Some quantum Lie algebras and their enveloping Hopf algebras} \label{aquantliealganditsenvalg}

We follow \cite{maj-gom2003} for the definitions and related notions that we are going to use in what follows. \newline

In this section we define a few particular examples of quantum Lie algebras, $\mathfrak{g}$, (in the sense of \cite{maj-gom2003}), for which, as we will show, their universal enveloping algebras, $U(\mathfrak{g})$, turn out to be the anyonic Hopf algebras that were defined in Section  \ref{anyonhpfalgwtheproppsi2=id}. 

\begin{example} \label{nonabelianquantliealginz2m2twogenerator}
Let $\mathcal{C}$ be the category of left $\mathbb{C} \mathbb{Z}_{n}$-modules, where $n=2m^2$, for some odd integer $m>1$. We define an object, $\mathfrak{g}=\bigoplus_{i=0}^{n-1} \mathfrak{g}_i$, in $\mathcal{C}$, by setting  $\mathfrak{g}_{im}= \mathbb{C}x_i$ and $\mathfrak{g}_{2im}= \mathbb{C}a_i$, for $i=1, 3, 5, 7, \cdots, 2m-1$, and $\mathfrak{g}_i= 0$, otherwise. We define a Lie bracket on $\mathfrak{g}$ by $[x_i,\, x_i]=a_i$, for $i=1, 3, 5, 7, \cdots, 2m-1$, and $[\,,\,]=0$ otherwise. Then, $(\mathfrak{g}, \psi_{\mathfrak{g} \otimes \mathfrak{g}}, [\, , \,])$ is a left quantum Lie algebra, in the sense of Definition 2.1 in \cite{maj-gom2003}.  
\end{example}

\begin{proof}
Since $\psi_{\mathfrak{g} \otimes \mathfrak{g}}(v\otimes w)=(-1)^{\frac{|v||w|}{m^2}} w\otimes v$, properties (1) and (3) in Definition 2.1 in \cite{maj-gom2003} are automatically satisfied. Also, since the only elements in $ker(id - \psi)$ are $a_i\otimes a_i$, and since by our definition $[a_i, a_i]=0$, property (4) is satisfied too. All one needs to check is property (2), the quantum Jacobi identity. On one hand, if one of the components in any bracket $[*\, ,* \,]$ is $a_i$, then $[*\, ,* \,]=0$. On the other hand, $[*\, ,* \,]$ is equal to either zero or $a_i$. These two statements imply that any bracket of the form $[*\, ,[* \, ,*] \,]$ or $[\,[*\, ,*] \, ,* \,]$ (appearing in (2)) is equal to zero. Therefore, in (2), we always get both sides equal to zero, and we are done.
\end{proof}

By Lemma 2.6 in \cite{maj-gom2003}, for the quantum Lie algebra, $\mathfrak{g}$, in the above example, the universal enveloping algebra, $U(\mathfrak{g})$, is a  Hopf algebra \cite{maja94-algandhopfalginbrddcats, maj95-book} with the standard structure. Let us see what  this Hopf algebra is. 
By the definition of $U(\mathfrak{g})$ \cite{maj-gom2003} and Formula \ref{siform2almostsuper}, we have, 
\[U(\mathfrak{g})=\frac{T(\mathfrak{g})}{(v\otimes w -(-1)^{\frac{|v||w|}{m^2}} w \otimes v-[v,\, w])}. \label{ugdefbyformula2.1} \]

This means that in $U(\mathfrak{g})$, $[v,w]=vw-wv$, when $|v||w|$ is an even multiple of $m^2$, and $[v,w]=vw+wv$, when $|v||w|$ is an odd multiple of $m^2$. When we apply these to the basis $x_i$, and $a_i$, we get, for $i,j=1, 3, 5, 7, \cdots, 2m-1$,
\[ a_i=[x_i, x_i]=2x_i^2, \qquad  0=[x_i, x_j]=x_ix_j+x_jx_i, \, i\neq j . \nonumber \]
Therefore, $U(\mathfrak{g})$ for the $\mathfrak{g}$ of Example \ref{nonabelianquantliealginz2m2twogenerator}, is, in fact, the anyonic Hopf algebra of Example \ref{anyhopfalginz2m2lessgeneral}. \newline

A special case of the above example is when we let $m=3$, as follows.

\begin{example} \label{nonabelianquantliealginz18twogenerator}
We define $\mathfrak{g}=\bigoplus_{i=0}^{17} \mathfrak{g}_i$, in the category $\mathcal{C}$ of left $\mathbb{C} \mathbb{Z}_{18}$-modules, by setting $\mathfrak{g}_3= \mathbb{C}x$, $\mathfrak{g}_6= \mathbb{C}a$, $\mathfrak{g}_9= \mathbb{C}y$, $\mathfrak{g}_{18}=\mathfrak{g}_0= \mathbb{C}b$ (remember $18 \equiv 0, \quad mod 18$), $\mathfrak{g}_{15}= \mathbb{C}z$, $\mathfrak{g}_{30}=\mathfrak{g}_{12}= \mathbb{C}c$ ($30 \equiv 12, \quad mod 18$), and $\mathfrak{g}_i= 0$, otherwise. We define a Lie bracket on $\mathfrak{g}$ by $[x, x]=a$, $[y, y]=b$, $[z, z]=c$, and $[\,,\,]=0$, otherwise. Then, $(\mathfrak{g}, \psi_{\mathfrak{g} \otimes \mathfrak{g}}, [\, , \,])$ is a left quantum Lie algebra, in the sense of Definition 2.1 in \cite{maj-gom2003}. Also, $U(\mathfrak{g})$ is the anyonic Hopf algebra of Example \ref{anyhopfalginz18lessgeneral}. 
\end{example}

Now, to cover the remaining Hopf algebras in Section \ref{anyonhpfalgwtheproppsi2=id} (except Example \ref{anyhopfalginz2m2moregeneral}), we present the following examples of abelian quantum Lie algebras. 

\begin{example} \label{quantliealginz2m2}
Let $\mathcal{C}$ be the category of left $\mathbb{C} \mathbb{Z}_{n}$-modules, where $n=2m^2$, for some odd integer $m>1$. We define $\mathfrak{g}=\bigoplus_{i=0}^{n} \mathfrak{g}_i$, in $\mathcal{C}$, by setting, $\mathfrak{g}_{im}= \mathbb{C}x_i$, for $i=1, 3, 5, 7, \cdots, 2m-1$, and $\mathfrak{g}_i= 0$, otherwise. We define the Lie bracket on $\mathfrak{g}$, to be zero. Then, it is very easy to see that $(\mathfrak{g}, \psi_{\mathfrak{g} \otimes \mathfrak{g}}, [\, , \,])$ is a quantum abelian Lie algebra.
\end{example}
For $\mathfrak{g}$ of this example, in a similar manner as before, and with easier calculations (as $\mathfrak{g}$ is abelian), we will see that $U(\mathfrak{g})$ is the anyonic Hopf algebra of Example \ref{anyhopfalginz2m2}.\newline

When we let $m=3$ we have, 
 
\begin{example} \label{quantliealginz18}
Let $\mathfrak{g}=\bigoplus_{i=0}^{17} \mathfrak{g}_i$, as an object in $\mathcal{C}$, the category $\mathcal{C}$ of left $\mathbb{C} \mathbb{Z}_{18}$-modules, be defined by setting $\mathfrak{g}_i= \mathbb{C}x_i$, for $i=1, 3,5$ (note, $|x_i|=3i$), and $\mathfrak{g}_i= 0$, otherwise. We define the Lie bracket on $\mathfrak{g}$, to be zero. The enveloping algebra $U(\mathfrak{g})$, for this quantum abelian Lie algebra, is the anyonic Hopf algebra of Example \ref{anyhopfalginz18}.
\end{example}

In Section \ref{hccforouranyonhofalg} we will use the results of this section, the results from \cite{pl4}, and an idea similar to the one used for super Hopf algebras in \cite{pl1hha}, to compute the braided periodic Hopf cyclic cohomology of $H=U(\mathfrak{g})$ from the above examples.

\section{Computing Hopf cyclic cohomology of some anyonic Hopf algebras} \label{hccforouranyonhofalg} 
For a quick review of the subject of (braided) Hopf cyclic cohomology we refer the reader to \cite{pl3}.\newline 

In \cite{pl4}, we showed that the main theorems of braided Hopf cyclic cohomology in \cite{pl1hha}, which were stated for symmetric monoidal categories, are still valid in non-symmetric categories, as long as for the main player, namely the braided Hopf algebra, $H$, we have $\psi^2_{H \otimes H}=id_{H \otimes H}$. In this section we want to use this fact to compute the (braided) periodic Hopf cyclic cohomology of the anyonic Hopf algebras from Section \ref{anyonhpfalgwtheproppsi2=id}.
In Section \ref{aquantliealganditsenvalg}, we proved that for any of those Hopf algebras, $H$ (except for the one in Example \ref{anyhopfalginz2m2moregeneral}), $H=U(\mathfrak{g})$, for some quantum Lie algebra, $\mathfrak{g}$. 

This is why the main result of this section is the following lemma. It is an analogue of Theorem 6.3 in \cite{pl1hha}, adapted for quantum Lie algebras in the category of anyonic vector spaces, within the family of Example \ref{n=2m2}.  \newline

For the reason mentioned in Remark \ref{remarkonweaklybrddobjects}, it is very straightforward to see that the proof of the following lemma is similar to the one of Theorem 6.3, given in detail in \cite{pl1hha}. All one needs to do is to replace all the degrees of the homogeneous elements in the proof of that theorem, such that a given degree, $|x|$, becomes $\frac{|x|}{m}$.  

\begin{lemma} \label{anyonverofHCCug}
Let $\mathcal{C}$ be the category of left $\mathbb{C} \mathbb{Z}_{n}$-modules, where $n=2m^2$ for some integer $m>1$. Let $(\mathfrak{g}, \psi_{\mathfrak{g} \otimes \mathfrak{g}}, [\, , \, ])$ be a quantum Lie algebra in $\mathcal{C}$, focused only in degrees, $km$, for integers $k\geq 0$ (i.e., within the family of Example \ref{n=2m2}). Then, analogous to the non-super case \cite{cm2,cm3,cm4,crn} and the super case \cite{pl1hha}, the relation,
\[ HP_{(\delta,1)}^*(U(\mathfrak{g}))=\bigoplus_{i=*(mod~2)} H_i (\mathfrak{g} ; \mathbb{C}_ \delta), \quad *=0,\,1 , \label{perhcsquantenvalg} \] 
holds, where $\delta$ is a character for $\mathfrak{g}$. Here, $HP_{(\delta,1)}^*(U(\mathfrak{g}))$ is
the periodic Hopf cyclic cohomology of the anyonic Hopf algebra, $H=U(\mathfrak{g})$, and $H_i (\mathfrak{g} ; \mathbb{C}_\delta) $ 
is the Lie algebra homology of $\mathfrak{g}$, with coefficients in the $\mathfrak{g}$-module $\mathbb{C}_\delta$. 
\end{lemma} 

The Lie algebra homologies, $H_i (\mathfrak{g} ; \mathbb{C}_\delta)$, in the above lemma, are computed by the following complex, which is an analogue of the Chevalley-Eilenberg complex, adapted for $\mathfrak{g}$ within the family of Example \ref{n=2m2}. 
\begin{displaymath}  \nonumber
\xymatrix{ \bigwedge^0\mathfrak{g}  & \ar[l]_{\delta}  \bigwedge^1\mathfrak{g}  & \ar[l]_{d}  \bigwedge^2\mathfrak{g} & \ar[l]_{d}  \bigwedge^3\mathfrak{g} & \ar[l]_{d}... }
\end{displaymath}
\begin{eqnarray} \label{CEcplxforanyonliealg}
d(x_1\wedge \cdots \wedge x_i) &=& \sum_{k=1}^i (-1)^{k+1+ \alpha_k} \, \delta(x_k)(x_1 \wedge...\wedge \widehat{x}_k \wedge...\wedge x_i) \\ \nonumber
&+&\sum_{1 \leq k<l \leq i} (-1)^{k+l+ \beta_{k,l}} \, ([x_k, x_l] \wedge x_1 \wedge...\wedge \widehat{x}_k \wedge...\wedge \widehat{x}_l \wedge...\wedge x_i) 
\end{eqnarray}
where, $\alpha_1:=0$, $\alpha_k := \frac{|x_k|(|x_1|+...+|x_{k-1}|)}{m^2}$ for $k>1$,  $\beta_{k,l} :=\alpha_k + \alpha_l - \frac{|x_k||x_l|}{m^2}$ for $k<l$,  and $\widehat{x}_k$ means omitted. Also, 
\[\bigwedge(\mathfrak{g})=\bigoplus_{i=0}^{\infty} \bigwedge^i \mathfrak{g} :=\frac{T(\mathfrak{g})}{( a\otimes b +(-1)^{\frac{|a||b|}{m^2}} b \otimes a)}. \nonumber\] \newline

Now, one can apply this lemma to any of the anyonic Hopf algebras in Section \ref{anyonhpfalgwtheproppsi2=id}. \newline

As an example, we close this paper by applying Lemma \ref{anyonverofHCCug}, in detail, to the anyonic Hopf algebra of Example \ref{anyhopfalginz18}, proved to be the enveloping algebra of the quantum Lie algebra of Example \ref{quantliealginz18}. \newline

In Formula \eqref{CEcplxforanyonliealg}, since $\delta=\varepsilon$, and the Lie algebra brackets are zero, all the boundary maps are zero. Thus $H_i (\mathfrak{g} ,\, \mathbb{C}_\varepsilon)=\bigwedge^i \mathfrak{g}$, where,  $\bigwedge^i \mathfrak{g}$, as a vector space, is generated by all elements of the form $\wedge^k x \wedge \wedge^l y \wedge \wedge^p z$ such that $k+l+p=i$. Here, by $\wedge^k x$ we mean $x \wedge x \wedge x \cdots \wedge x$, $k$ times. 

On the other hand, since in $\mathfrak{g}=\bigoplus_{i=0}^{17} \mathfrak{g}_i$, $|x|=3$, $|y|=9$, and $|z|=15$, we have $|\wedge^k x \wedge \wedge^l y \wedge \wedge^p z|=3k+9l+15p=3(i+(2l+4p))$,  (mod 18).
This implies that, as an anyonic vector space, $H_i (\mathfrak{g})= \bigoplus_{t=0}^{17} H_i (\mathfrak{g})_t$ is focused only in degrees that are even multiples of $3$ when $i$ is even, and in degrees that are odd multiples of $3$ when $i$ is odd.

Therefore by Formula \eqref{perhcsquantenvalg}, as an anyonic vector space,  
$HP_{(\varepsilon,1)}^0(U(\mathfrak{g}))=\bigoplus_{i=0(mod~2)} \bigwedge^i \mathfrak{g}$ is focused in degrees that are even multiples of $3$, with infinitely many elements in the basis of each degree. Similarly, 
$HP_{(\varepsilon,1)}^1(U(\mathfrak{g}))=\bigoplus_{i=1(mod~2)} \bigwedge^i \mathfrak{g}$ is focused in degrees that are odd multiples of $3$, with infinitely many elements in the basis of each degree. 


\end{document}